\documentclass[review]{elsarticle}

\usepackage{lineno}
\usepackage[hypertexnames=false]{hyperref}
\modulolinenumbers[5]
\usepackage{graphicx} 
\usepackage{amsmath,amssymb}
\usepackage{subcaption}
\usepackage{algorithm}
\usepackage{algpseudocode}
\newtheorem{example}{Problem}[section]
\newtheorem{remark}{Remark}[section]
\newtheorem{theorem}{Theorem}[section]
\journal{Journal of Computational and Applied Mathematics.}
\usepackage{cleveref}

\begin{document}
\nolinenumbers
\begin{frontmatter}

\title{Non-stationary Anderson acceleration with optimized damping \tnoteref{mytitlenote}}
\tnotetext[mytitlenote]{\textbf{Funding:} This work was partially supported by the National Natural Science Foundation of China [grant number 12001287]; the Startup Foundation for Introducing Talent of Nanjing University of Information Science and Technology [grant number 2019r106]}


\author[mymainaddress,mysecondaryaddress]{Kewang Chen\corref{mycorrespondingauthor}}
\cortext[mycorrespondingauthor]{Corresponding author}
\ead{kwchen@nuist.edu.cn}

\author[mysecondaryaddress]{Cornelis Vuik}
\ead{c.vuik@tudelft.nl} \ead[url]{https://homepage.tudelft.nl/d2b4e/}

\address[mymainaddress]{College of Mathematics and Statistics, Nanjing University of Information Science and Technology, Nanjing, 210044, China.}
\address[mysecondaryaddress]{Delft Institute of Applied Mathematics, Delft University of Technology, Delft, 2628XE, the Netherlands.}

\begin{abstract}
Anderson acceleration (AA) has a long history of use and a strong recent interest due to its potential ability to dramatically improve the linear convergence of the fixed-point iteration. Most authors are simply using and analyzing the stationary version of Anderson acceleration (sAA) with a constant damping factor or without damping. Little attention has been paid to nonstationary algorithms. However, damping can be useful and is sometimes crucial for simulations in which the underlying fixed-point operator is not globally contractive. The role of this damping factor has not been fully understood. In the present work, we consider the non-stationary Anderson acceleration algorithm with optimized damping (AAoptD) in each iteration to further speed up linear and nonlinear iterations by applying one extra inexpensive optimization. We analyze this procedure and develop an efficient and inexpensive implementation scheme. We also show that, compared with the stationary Anderson acceleration with fixed window size $sAA(m)$, optimizing the damping factors is related to dynamically packaging $sAA(m)$ and $sAA(1)$ in each iteration (alternating window size $m$ is another direction of producing non-stationary AA). Moreover, we show by extensive numerical experiments that, in the case a larger window size is needed, the proposed non-stationary Anderson acceleration with optimized damping procedure often converges much faster than stationary AA with constant damping or without damping. When the window size is very small ($m\leq3$ was typically used, especially in the early days of application), AAoptD and AA are comparable. Lastly, we observed that when the system is overdamped (i.e. the damping factor is close to the lower bound zero), inconsistency may occur. So there is some trade-off between stability and speed of convergence. We successfully solve this problem by further restricting damping factors bound away from zero.
\end{abstract}


\end{frontmatter}

\Crefname{figure}{Fig.}{Figs.}
\section{Introduction}
\label{sec:intro}
In this part, we first give a literature review on Anderson Acceleration method. Then we discuss our main motivations and the structure for the present paper. To begin with, let us consider the nonlinear acceleration for the following general fixed-point problem 
$$x=g(x),\  g : R^n\rightarrow R^n$$
or its related nonlinear equations problem 
 $$f(x)=x-g(x)=0.$$
The associated basical fixed-point iteration is given in \Cref{alg:picard}.
\begin{algorithm}
\caption{Picard iteration}
\label{alg:picard}
\begin{algorithmic}
\State{Given: $x_0$.}
\For{$k=0,1,2,\cdots $}
\State{Set $x_{k+1}=g(x_k)$.}
\EndFor
\end{algorithmic}
\end{algorithm}

The main concern related to this basic fixed-point iteration is that the iterates may not converge or may converge extremely slowly (only linear convergent). Therefore, various acceleration methods are proposed to alleviate this slow convergence problem. Among these algorithms, one popular acceleration procedure is called the Anderson acceleration method  \cite{An1965}. For the above basic Picard iteration, the usual general form of Anderson acceleration with damping is given in \Cref{alg:Anderson}.
\begin{algorithm}
\caption{Anderson acceleration: $AA(m)$}
\label{alg:Anderson}
\begin{algorithmic}
\State{Given: $x_0$ and $m\geq 1$.}
\State{Set: $x_{1}=g(x_0).$}
\For{$k=0,1,2,\cdots $}
\State{Set: $m_{k}=\min\{m,k\}$.}
\State{Set: $F_k=(f_{k-m_k},\cdots,f_k)$, where $f_i=g(x_i)-x_i$.}
\State{Determine: $\alpha^{(k)}=\left(\alpha_0^{(k)},\cdots,\alpha_{m_k}^{(k)}\right)^{T}$ that solves }
\State{$\ \ \ \ \ \ \ \ \ \ \ \displaystyle \min_{\alpha=(\alpha_0,\cdots,\alpha_{m_k})^{T}}\|F_k\alpha\|_2$ $\ s.\ t.$ $\displaystyle\sum_{i=0}^{m_k}\alpha_i=1.$}
\State{Set: $\displaystyle x_{k+1}=(1-\beta_k)\sum_{i=0}^{m_k}\alpha_{i}^{(k)}x_{k-m_k+i}+\beta_k\sum_{i=0}^{m_k}\alpha_{i}^{(k)}g(x_{k-m_k+i})$.}
\EndFor
\end{algorithmic}
\end{algorithm}
In the above algorithm, $f_k$ is the residual for the $k$th iteration; $m$ is the window size which indicates how many history residuals will be used in the algorithm. The value of $m$ is typically no larger than $3$ in the early days of applications and now this value could be as large as up to 100, see \cite{An2019}. It is usually a fixed number during the procedure, varying $m$ can also make the algorithm to be non-stationary. We will come back to this point in section \Cref{sec:main}; $\beta_k\in (0,1]$ is a damping factor (or a relaxation parameter) at $k$th iteration. We have, for a fixed window size $m$:
\begin{equation*}
   \beta_k=
    \begin{cases}
      1,  & \text{no damping,}\\
      \beta, \  (\text{a constant independent of $k$}) & \text{stationary AA,}\\
      \beta_k,\ (\text{depending on $k$}) \  & \text{non-stationary AA.}
    \end{cases}       
\end{equation*}
The constrained optimization problem can also be formulated as an equivalent unconstrained least-squares problem \cite{ToKe2015,wa2011}:
\begin{equation}
\displaystyle \min_{(\omega_1,\cdots,\omega_{m_k})^{T}}\left\|f_k+\sum_{i=1}^{m_k}\omega_i(f_{k-i}-f_k)\right\|_2
\end{equation}
One can easily recover the original problem by setting
$$\omega_0=1-\sum_{i=1}^{m_k}\omega_i.$$
This formulation of the linear least-squares problem is not optimal for implementation, we will discuss this in more detail in \Cref{sec:proof}.

Anderson acceleration method dates back to the 1960s. In 1962, Anderson \cite{An1965} developed a technique for accelerating the convergence of the Picard iteration associated with a fixed-point problem which is called Extrapolation Algorithm. This technique is now called Anderson Acceleration (AA) in the applied mathematics community and Anderson Mixing in the physics and chemistry communities. This method is ``essentially'' (or nearly) similar to the nonlinear GMRES method or Krylov acceleration \cite{CaMi1998,Mi2005,OoWa2000,WaOo1997} and the direct inversion on the iterative subspace method (DIIS) \cite{LiYa2013,Pu1980,Pu1982}. And it is also in a broad category with methods based on quasi-Newton updating \cite{EiNe1987,Ey1996,FaSa2009,Ha2010,Ya2009}. However, unlike Newton-like methods, AA does not require the computation or approximation of Jacobians or Jacobian-vector products which could be an advantage.

Although the Anderson acceleration method has been around for decades, convergence analysis has been reported in the literature only recently. Fang and Saad \cite{FaSa2009} had clarified a remarkable relationship of AA to quasi-Newton methods and extended it to define a broader Anderson family method. Later, Walker and Ni \cite{WaNi2011} showed that, on linear problems, AA without truncation is ``essentially equivalent'' in a certain sense to the GMRES method. For the linear case, Toth and Kelley \cite{ToKe2015} first proved the stationary version of AA (sAA) without damping is locally r-linearly convergent if the fixed point map is a contraction and the coefficients in the linear combination remain bounded. This work was later extended by Evens et al. \cite{Ev2020} to AA with damping and the authors proved the new convergence rate is $\theta_k((1-\beta_{k-1})+\beta_{k-1}\kappa)$, where $\kappa$ is the Lipschitz constant for the function $g(x)$ and $\theta_k$ is the ratio quantifying the convergence gain provided by AA in step $k$. However, it is not clear how $\theta_k$ may be evaluated or bounded in practice and how it may translate to improved asymptotic convergence behavior in general. In 2019, Pollock et al. \cite{Po2019} applied sAA to the Picard iteration for solving steady incompressible Navier–Stokes equations (NSE) and proved that the acceleration improves the convergence rate of the Picard iteration. Then, De Sterck \cite{DeHe2021} extended the result to more general fixed-point iteration $x=g(x)$, given knowledge of the spectrum of $g'(x)$ at fixed-point $x^*$ and Wang et al. \cite{WaSt2021} extended the result to study the asymptotic linear convergence speed of sAA applied to Alternating Direction Method of Multipliers (ADMM) method. Sharper local convergence results of AA remain a hot research topic in this area. More recently, Zhang et al. \cite{Zh2020} proved a global convergent result of type-{I} Anderson acceleration for nonsmooth fixed-point iterations without resorting to line search or any further assumptions other than nonexpansiveness. For more related results about Anderson acceleration and its applications, we refer the interested readers to \cite{An2019,BiKe2021,Br2015,YuLi2018,To2017,WeGa2019,Ya2021} and references therein.

As mentioned above, the local convergence rate $\theta_k((1-\beta_{k-1})+\beta_{k-1}\kappa)$ at stage $k$ is closely related to the damping factor $\beta_{k-1}$. However, questions like how to choose those damping values in each iteration \cite{An2019} and how it will affect the global convergence of the algorithm have not been deeply studied. Besides, AA is often combined with globalization methods to safeguard against erratic convergence away from a fixed point by using damping. One similar idea in the optimization context for nonlinear GMRES is to use line search strategies \cite{De2012}. This is an important strategy but not yet fully explored in the literature. Moreover, the early days of Anderson Mixing method (the 1980s, for electronic structure calculations) initially dictated the window size $m\leq 3$ due to the storage limitations and costly $g$ evaluations involving large $N$. However, in recent years and a broad range of contexts, the window size $m$ ranging from $20$ to $100$ has also been considered by many authors. For example, Walker and Ni \cite{WaNi2011} used $m=50$ in solving the nonlinear Bratu problem. A natural question will be should we try to further steep up Anderson acceleration method or try to use a larger size of the window? No such comparison results have been reported. Motivated by the above works, in this paper, we propose, analyze and numerically study non-stationary Anderson acceleration with optimized damping to solve fixed-point problems. The goal of this paper is to explore the role of damping factors in non-stationary Anderson acceleration.

The paper is organized as follows. Our new algorithms and analysis are in
\Cref{sec:main}, the implementation of the new algorithm is in \Cref{sec:proof}, experimental
results and discussion are in \Cref{sec:experiments}. Conclusions follow in
\Cref{sec:conclusions}.

\section{Anderson acceleration with optimized dampings}
\label{sec:main}

In this section, we focus on developing the algorithm for Anderson acceleration with optimized dampings at each iteration and studying its convergence rate explicitly.
\begin{eqnarray}\label{eq:1}
x_{k+1}&=&(1-\beta_k)\sum_{i=0}^{m_k}\alpha_{i}^{(k)}x_{k-m_k+i}+\beta_k\sum_{i=0}^{m_k}\alpha_{i}^{(k)}g(x_{k-m_k+i})\nonumber\\
&=&\sum_{i=0}^{m_k}\alpha_{i}^{(k)}x_{k-m_k+i}+\beta_k\left(\sum_{i=0}^{m_k}\alpha_{i}^{(k)}g(x_{k-m_k+i})-\sum_{i=0}^{m_k}\alpha_{i}^{(k)}x_{k-m_k+i}\right).
\end{eqnarray}
Define the following averages given by the solution $\alpha^{k}$ to the optimization problem by
\begin{equation}\label{eq:2}
x_{k}^{\alpha}=\sum_{i=0}^{m_k}\alpha_{i}^{(k)}x_{k-m_k+i},\ \ \ \tilde{x}_{k}^{\alpha}=\sum_{i=0}^{m_k}\alpha_{i}^{(k)}g(x_{k-m_k+i}).
\end{equation}
Then \eqref{eq:1} becomes
\begin{equation}
x_{k+1}=x_{k}^{\alpha}+\beta_k(\tilde{x}_{k}^{\alpha}-x_{k}^{\alpha}).
\end{equation}
A natural way to choose ``best'' $\beta_k$ at this stage is that choosing $\beta_k$ such that $x_{k+1}$ gives a minimal residual. This is similar to the original idea of Anderson acceleration with window size equal to one. So we just need to solve the following unconstrained optimization problem:
\begin{equation}\label{eq:3}
\min_{\beta_k}\|x_{k+1}-g(x_{k+1})\|_2=\min_{\beta_k}\|x_{k}^{\alpha}+\beta_k(\tilde{x}_{k}^{\alpha}-x_{k}^{\alpha})-g(x_{k}^{\alpha}+\beta_k(\tilde{x}_{k}^{\alpha}-x_{k}^{\alpha}))\|_2.
\end{equation}
Noting the fact that
\begin{eqnarray}
g(x_{k}^{\alpha}+\beta_k(\tilde{x}_{k}^{\alpha}-x_{k}^{\alpha}))&\approx&g(x_{k}^{\alpha})+\beta_k\frac{\partial g}{\partial x}\Big|_{x_k^{\alpha}}(\tilde{x}_{k}^{\alpha}-x_{k}^{\alpha})\nonumber\\
&\approx&g(x_{k}^{\alpha})+\beta_k\left(g(\tilde{x}_{k}^{\alpha})-g(x_{k}^{\alpha})\right).
\end{eqnarray}
Therefore, \eqref{eq:3} becomes
\begin{eqnarray}\label{eq:4}
&&\min_{\beta_k}\|x_{k+1}-g(x_{k+1})\|_2\nonumber\\
&=&\min_{\beta_k}\|x_{k}^{\alpha}+\beta_k(\tilde{x}_{k}^{\alpha}-x_{k}^{\alpha})-g(x_{k}^{\alpha}+\beta_k(\tilde{x}_{k}^{\alpha}-x_{k}^{\alpha}))\|_2\nonumber\\
&\approx&\min_{\beta_k}\|x_{k}^{\alpha}+\beta_k(\tilde{x}_{k}^{\alpha}-x_{k}^{\alpha})-\left[g(x_{k}^{\alpha})+\beta_k(g(\tilde{x}_{k}^{\alpha})-g(x_{k}^{\alpha}))\right]\|_2\nonumber\\
&\approx&\min_{\beta_k}\|\left(x_{k}^{\alpha}-g(x_{k}^{\alpha})\right)-\beta_k\left[(g(\tilde{x}_{k}^{\alpha})-g(x_{k}^{\alpha}))-(\tilde{x}_{k}^{\alpha}-x_{k}^{\alpha})\right]\|_2.
\end{eqnarray}
Thus, we just need to calculate the projection
\begin{equation}
\beta_k=\Big|\frac{\left(x_{k}^{\alpha}-g(x_{k}^{\alpha})\right)\cdot\left[\left(x_{k}^{\alpha}-g(x_{k}^{\alpha})\right)-(\tilde{x}_{k}^{\alpha}-g(\tilde{x}_{k}^{\alpha}))\right]}{\|\left[\left(x_{k}^{\alpha}-g(x_{k}^{\alpha})\right)-(\tilde{x}_{k}^{\alpha}-g(\tilde{x}_{k}^{\alpha}))\right]\|_2}\Big|.
\end{equation}
Set
$${r_p}=\left(x_{k}^{\alpha}-g(x_{k}^{\alpha})\right),\ \ {r_q}=\left(\tilde{x}_{k}^{\alpha}-g(\tilde{x}_{k}^{\alpha})\right),$$
we have
\begin{equation}
\beta_k=\left|\frac{(r_p-r_q)^{T}r_p}{\|r_p-r_q\|_2}\right|.
\end{equation}
We will discuss how much work is needed to calculate this $\beta_k$ in \Cref{sec:proof}. 
Finally, our analysis leads to the following non-stationary Anderson acceleration algorithm with optimized damping: $AAoptD(m)$.
\begin{algorithm}
\caption{Anderson acceleration with optimized dampings: $AAoptD(m)$}
\label{alg:AAoptD}
\begin{algorithmic}
\State{Given: $x_0$ and $m\geq 1$.}
\State{Set: $x_{1}=g(x_0).$}
\For{$k=0,1,2,\cdots $}
\State{Set: $m_{k}=\min\{m,k\}$.}
\State{Set: $F_k=(f_{k-m_k},\cdots,f_k)$, where $f_i=g(x_i)-x_i$.}
  \State{Determine: $\alpha^{(k)}=\left(\alpha_0^{(k},\cdots,\alpha_{m_k}^{(k)}\right)^{T}$ that solves }
\State{$\ \ \ \ \ \ \ \ \ \ \ \displaystyle \min_{\alpha=(\alpha_0,\cdots,\alpha_{m_k})^{T}}\|F_k\alpha\|_2$ $\ s.\ t.$ $\displaystyle\sum_{i=0}^{m_k}\alpha_i=1.$}
\State{Set: $\ \ \displaystyle x_{k}^{\alpha}=\sum_{i=0}^{m_k}\alpha_{i}^{(k)}x_{k-m_k+i},\ \ \ \tilde{x}_{k}^{\alpha}=\sum_{i=0}^{m_k}\alpha_{i}^{(k)}g(x_{k-m_k+i}).$}
\State{Set: $\ \ \displaystyle {r_p}=\left(x_{k}^{\alpha}-g(x_{k}^{\alpha})\right),\ \ {r_q}=\left(\tilde{x}_{k}^{\alpha}-g(\tilde{x}_{k}^{\alpha})\right)$.}
\State{Set: $\ \ \ \ \ \ \ \ \ \ \ \displaystyle \beta_k=\frac{(r_p-r_q)^{T}r_p}{\|r_p-r_q\|_2}$.}
\State{Set: $ \ \ \displaystyle x_{k+1}=x_{k}^{\alpha}+\beta_k(\tilde{x}_{k}^{\alpha}-x_{k}^{\alpha})$.}
\EndFor
\end{algorithmic}
\end{algorithm}

\begin{remark}
As mentioned in \Cref{sec:intro}, changing the window size $m$ at each iteration can also make a stationary Anderson acceleration to be non-stationary. Comparing with the stationary Anderson acceleration with fixed window $sAA(m)$, our proposed nonstationary procedure ($AAoptD(m)$) of choosing optimal $\beta_k$ is somewhat related to packaging $sAA(m)$ and $sAA(1)$ in each iteration in a cheap way. Combining $sAA(m)$ with $sAA(1)$ can provide really good outcomes, especially in the case when larger $m$ is needed. We will discuss this in detail for the numerical results in \Cref{sec:experiments}.
\end{remark}

\begin{remark}
Here this optimized damping step is a ``local optimal'' strategy at $k$th iteration. It usually will speed up the convergence rate compared with an undamped one, but not always. Because in $(k+1)$th iteration, it uses a combination of all previous m history information. Moreover, when $\beta_k$ is very close to zero, the system is over-damped, which, sometimes, may also slow down the convergence speed. We may need to further modify our $\beta_k$. See more discussion in our numerical results in \Cref{sec:experiments}.
\end{remark}


Lastly, we summarize the convergence results with damping in \Cref{thm:bigthm}. The proof of this theorem can be found in \cite{Ev2020}. 

\begin{theorem}\cite{Ev2020} \label{thm:bigthm}Assume that $g: R^n\rightarrow R^n$ is uniformly Lipschitz continuously differentiable and there exists $\kappa\in(0,1)$ such that $\|g(y)-g(x)\|_2\leq\kappa\|y-x\|_2$ for all $x, y\in R^n$. Suppose also that $\exists M$ and $\epsilon>0$ such that for all $k>m$, $\sum_{i=0}^{m-1}|\alpha_i|<M$ and $|\alpha_m|\geq\epsilon$. Then
\begin{equation}
\|f(x_{k+1})\|_2\leq\theta_{k+1}\left[(1-\beta_k)+\kappa\beta_k\right]\|f(x_k)\|_2+\sum_{i=0}^{m} O(\|f(x_{k-m+i})\|_2^2),
\end{equation}
where
$$\theta_{k+1}=\frac{\|\sum_{i=0}^{m}\alpha_if(x_{k-m+i})\|_2}{\|f(x_k)\|_2}.$$
\end{theorem}

\section{Implementation}\label{sec:proof}
For implementation, we mainly follow the path in \cite{wa2011} and modify it as needed. We first briefly review the implementation of AA without damping. Then we focus on how to implement the optimized damping problem efficiently and accurately.

The constrained linear least-squares problem in \Cref{alg:Anderson} can be solved in many ways. Here we rewrite it into an equivalent unconstrained form which can be solved efficiently by using QR factorizations. We define $\Delta f_i=f_{i+1}-f_i$ for each $i$ and set $\mathcal{F}_k=(\Delta f_{k-m_k},\cdots,\Delta f_{k-1})$, then the least-squares problem is equivalent to 
$$\min_{\gamma=(\gamma_0,\cdots,\gamma_{m_{k-1}})^T}\|f_k-\mathcal{F}_k\gamma\|_2,$$
where $\alpha$ and $\gamma$ are related by $\alpha_0=\gamma_0, \alpha_i=\gamma_i-\gamma_{i-1}$ for $1\leq i\leq m_k-1$, and $\alpha_{m_k}=1-\gamma_{m_k-1}.$ We assume $\mathcal{F}$ has a thin $QR$ decomposition i.e.,  $\mathcal{F}_k=Q_kR_k$ with $Q_k\in \mathcal{R}^{n\times m_k}$ and $R_k\in \mathcal{R}^{m_k\times m_k}$, for which the solution of the least-squares problem is obtained by solving the $m_k\times m_k$ triangular system $R_k\gamma=Q_k^{T}f_k$. As the algorithm proceeds, the successive least-squares problems can be solved efficiently by updating the factors in the decomposition. 

Assume that $\gamma^{k}=(\gamma_0^{k},\cdots,\gamma_{m_k-1}^{k})^{T}$ is the solution to the above modified form of Anderson acceleration, we have
$$x_{k+1}=g(x_k)-\sum_{i=0}^{m_k-1}\gamma_i^{k}\left[g(x_{k-m_k+i+1})-g(x_{k-m_k+i})\right]=g(x_k)-\mathcal{G}_k\gamma^k,$$
where $\mathcal{G}_k=(\Delta g_{g_{k-m_k}},\cdots,\Delta g_{k-1})$ with $\Delta g_i=g(x_{i+1}-g(x_i))$ for each $i$. For Anderson acceleration with damping
\begin{eqnarray*}
x_{k+1}&=&(1-\beta_k)\sum_{i=0}^{m_k}\alpha_{i}^{(k)}x_{k-m_k+i}+\beta_k\sum_{i=0}^{m_k}\alpha_{i}^{(k)}g(x_{k-m_k+i})\nonumber\\
&=&\sum_{i=0}^{m_k}\alpha_{i}^{(k)}x_{k-m_k+i}+\beta_k\left(\sum_{i=0}^{m_k}\alpha_{i}^{(k)}g(x_{k-m_k+i})-\sum_{i=0}^{m_k}\alpha_{i}^{(k)}x_{k-m_k+i}\right).
\end{eqnarray*}
Follow the idea in \cite{wa2011}, we have
\begin{equation}\label{eq:kw1}
\sum_{i=0}^{m_k}\alpha_{i}^{(k)}g(x_{k-m_k+i})=g(x_k)-\mathcal{G}_k\gamma^{k},
\end{equation}
\begin{equation}\label{eq:kw2}
\sum_{i=0}^{m_k}\alpha_{i}^{(k)}x_{k-m_k+i}=\left(g(x_k)-\mathcal{G}_k\gamma^{k}\right)-\left(f_k-\mathcal{F}_k\gamma^k\right).
\end{equation}
Then this can be achieved equivalently using the following strategy:

Step 1: Compute the undamped iterate $x_{k+1}=g(x_k)-\mathcal{G}_k\gamma^{k}$.

Step 2: Update $x_{k+1}$ again by
$$x_{k+1}\leftarrow x_{k+1}-(1-\beta_k)\left(f_k-QR\gamma^k\right).$$
Now we talk about how to efficiently calculate $\beta_k$ as described in \Cref{alg:AAoptD}. Taking benefit of the QR decomposition in the first optimization problem and noting \eqref{eq:kw1} and \eqref{eq:kw2}, we have

$$\tilde{x}_{k}^{\alpha}=\sum_{i=0}^{m_k}\alpha_{i}^{(k)}g(x_{k-m_k+i})=g(x_k)-\mathcal{G}_k\gamma^{k},$$
$$x_{k}^{\alpha}=\sum_{i=0}^{m_k}\alpha_{i}^{(k)}x_{k-m_k+i}=\tilde{x}_{k}^{\alpha}-\left(f_k-\mathcal{F}_k\gamma^k\right).$$
Then we could calculate optimized $\beta_k$ by doing two extra function evaluations and two dot products, which are not very expensive:
$${r_p}=\left(x_{k}^{\alpha}-g(x_{k}^{\alpha})\right),\ \ {r_q}=\left(\tilde{x}_{k}^{\alpha}-g(\tilde{x}_{k}^{\alpha})\right),\ \
\beta_k=\left|\frac{(r_p-r_q)^{T}r_p}{\|r_p-r_q\|_2}\right|. $$
In practice, when $x_k$ is very close to the fixed-point $x^*$, scientific computing errors may arise in calculating these two high dimension vectors $r_p$ and $r_p-r_q$. Thus we normalize these two vectors first, then calculate $\beta_k$ by simply doing a dot product.

\section{Experimental results and discussion}
\label{sec:experiments} In this section, we numerically compare the performance of this non-stationary AAoptD with sAA (with constant damping or without damping). The first part contains examples where larger window sizes $m$ are needed in order to accelerate the iteration. The second part consists of some examples where small window sizes are working very well. All these experiments are done in MATLAB 2021b environment. MATLAB codes are available upon request to the authors.

This first example is from Walker and Ni's \cite{WaNi2011} paper, where a stationary Anderson acceleration with window size $m=50$ is used to solve the Bratu problem. This problem has a long history, we refer the reader to Glowinski et al. \cite{Gl1985} and Pernice and Walker \cite{Pe1998}, and the references in those papers. It is not a difficult problem for Newton-like solvers. 
\begin{example}\textbf{The\ Bratu\ problem.} The Bratu problem is a nonlinear PDE boundary value problem as follows:
\begin{eqnarray*}
\Delta u +\lambda\ e^{u}&=&0,\ \ in \ \ D=[0,1]\times[0,1],\nonumber\\
u&=&0,\ \ on \ \ \partial D. 
\end{eqnarray*}
\end{example}
In this experiment, we used a centered-difference discretization on a $32\times 32$, $64\times 64$ and $128\times 128$ grid, respectively. We take $\lambda=6$ in the Bratu problem and use the zero initial approximate solution in all cases. We also applied preconditioning such that the basic Picard iteration still works. The preconditioning matrix that we used here is the diagonal inverse of the matrix $A$, where $A$ is a matrix for the discrete Laplace operator.

\begin{figure}[htbp]
  \centering
  \includegraphics[height=3.4in, width=4.5in]{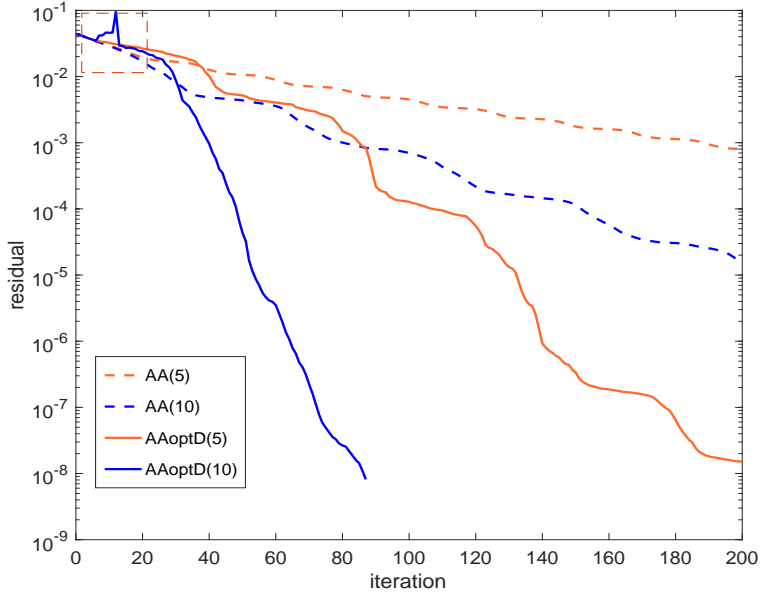}
  \caption{Compare AA and AoptD for solving nonlinear Bratu problems.}
  \label{fig:f1}
\end{figure}
\begin{figure}[htbp]
  \centering
  \includegraphics[height=1.8in, width=4.5in]{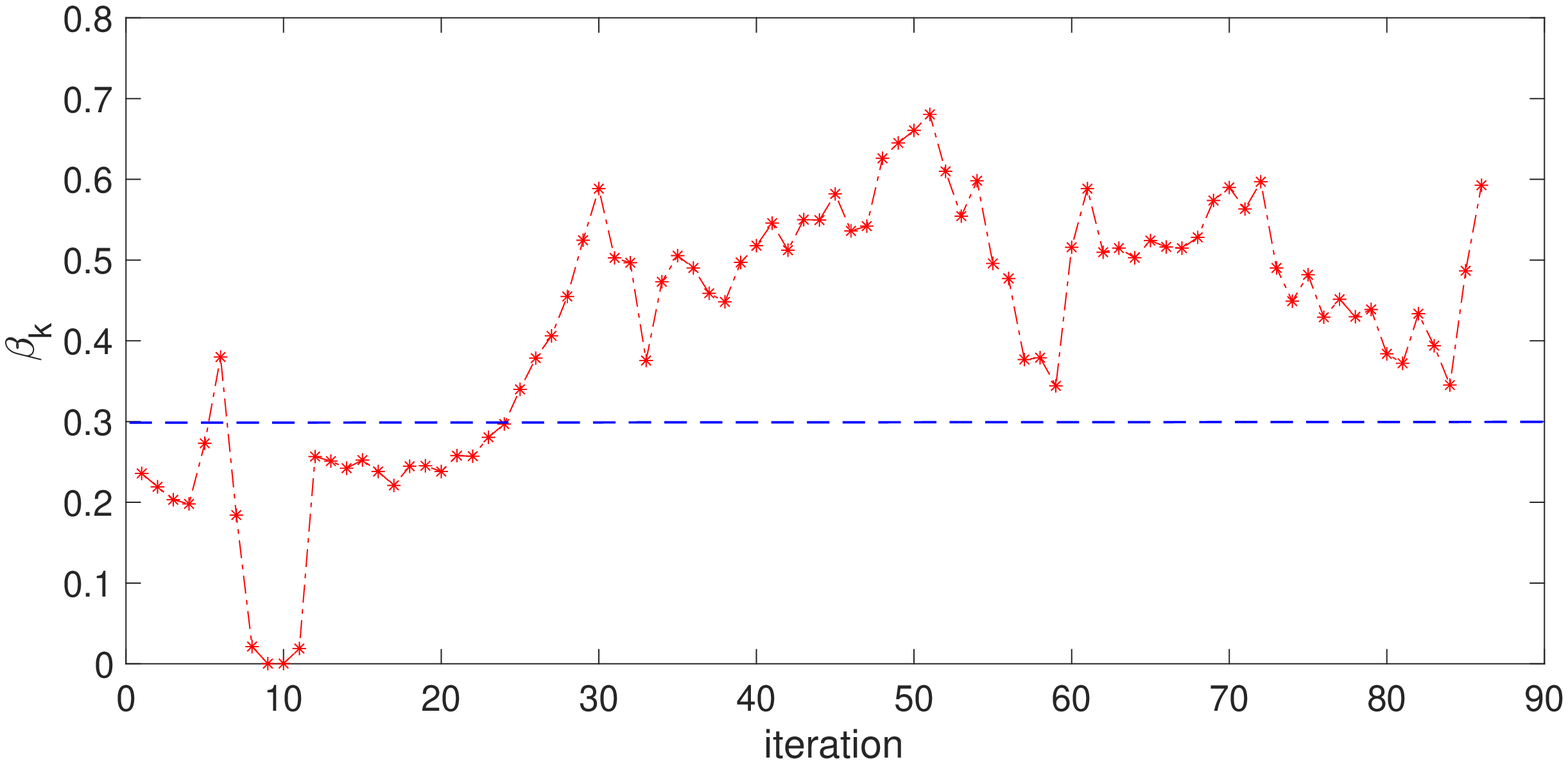}
  \caption{Optimal damping factors in each iteration for $m=10$.}
  \label{fig:f2}
\end{figure}

The results are shown in the following figures. In \Cref{fig:f1}, we plot the results of applying $AA(m)$ and $AAoptD(m)$ to accelerate Picard iteration with $m=5$ and $m=10$ on a grid of $32\times 32$. As we see from the picture, $AA(5)$ and $AA(10)$ does not accelerate the convergence speed very much. $AAoptD(5)$ and $AAoptD(10)$ perform much better than $AA(5)$ and $AA(10)$. However, we also notice that there are some inconsistencies and stagnations in $AAoptD(m)$. Thus we go further to plot the $\beta_k$ values that are used in each iteration, see \Cref{fig:f2}. From \Cref{fig:f2} we see that: for $AAoptD(10)$, some optimized damping factors are below $0.3$ (see the dashed line). As we know, the damping factor $\beta_k\in(0,1]$ and $\beta_k=1$ means no damping. Thus small $\beta_k$ may cause an over-damping phenomenon, which might be the reason for small inconsistencies observed in \Cref{fig:f1}; Similarly, we see that the residual of $AAoptD(10)$ in \Cref{fig:f1} is not decreasing consistently around $10$th iteration (see the read dashed square region in \Cref{fig:f1}), where the corresponding $\beta_k$ values are super close to zero as shown in \Cref{fig:f2}. 

\begin{figure}[htbp]
  \centering
  \includegraphics[height=1.8in, width=4.5in]{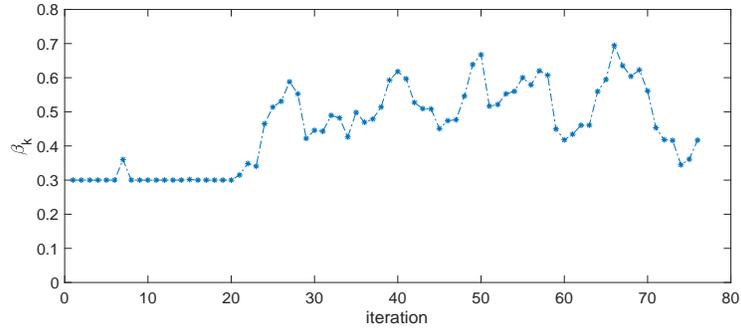}
  \caption{Modified optimal damping factors: $\hat{\beta}_k=\max\{\beta_k,\eta\}$ with $\eta=0.3$}
  \label{fig:new_f5}
\end{figure}

\begin{figure}[htbp]
  \centering
  \includegraphics[height=3.4in, width=4.5in]{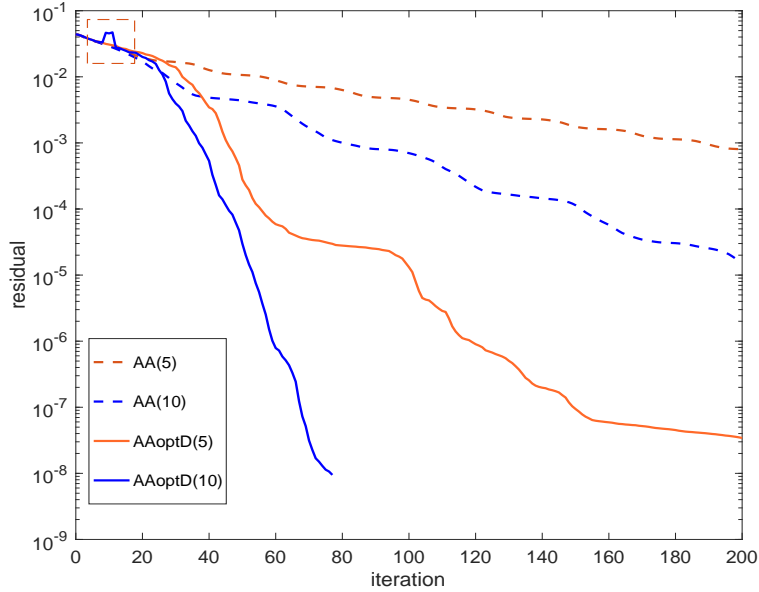}
  \caption{Solving nonlinear Bratu problems: $\hat{\beta}_k=\max\{\beta_k,\eta\}$ with $\eta=0.3$}
  \label{fig:new_f6}
\end{figure}

To balance the over-damping effect, we bound these $\beta_k$ away from zero. The first strategy we propose is to use
\begin{equation}\label{damp2}
\hat{\beta}_k =\max\{\beta_k,\eta\},
\end{equation}
where $\eta$ is a small positive number such that $0<\eta<0.5$. For example, to reduce the over-damping effect, we take $\eta=0.3$ in \eqref{damp2} as a lower bound. We plot the new $\beta_k$ values at each iteration in \Cref{fig:new_f5}. There are no $\beta_k$ values less than $0.3$ anymore. The corresponding results are in \Cref{fig:new_f6}. Compared with the results in \Cref{fig:f1}, we see that there is less stagnation (see the red dashed square region in \Cref{fig:new_f6}) and onvergence is also faster. We also note that the $\beta_k$ values in \Cref{fig:new_f5} differs a lot from the values of $\beta_k$ in \Cref{fig:f2}. Because changing $\beta_k$ in previous iterations will affect the later ones.

\begin{figure}[htbp]
  \centering
  \includegraphics[height=3.5in, width=4.5in]{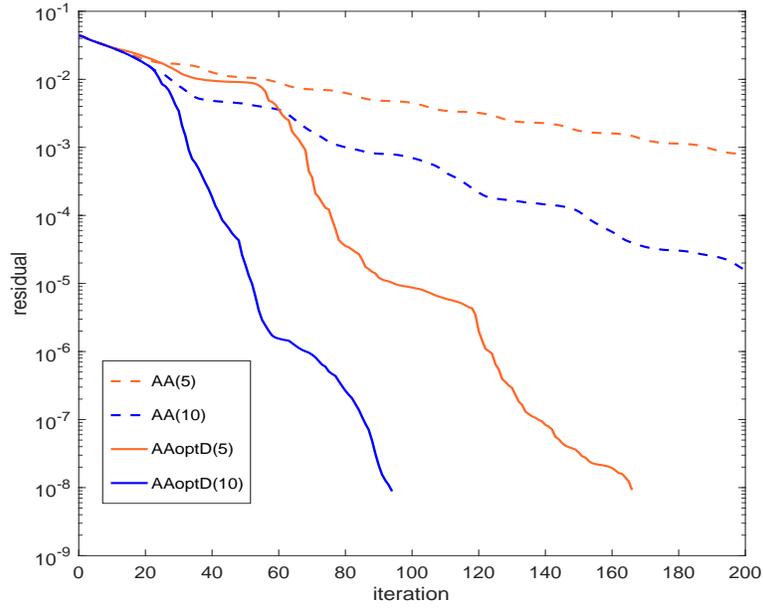}
  \caption{Solving nonlinear Bratu problems: $\hat{\beta}_k=1-\beta_k$ when $\beta_k<0.3$.}
  \label{fig:f4}
\end{figure}

Although the results in \Cref{fig:new_f6} are better than those in \Cref{fig:f1}, we notice that there are still some inconsistencies in the red dashed square region. To further smooth out these inconsistencies, we change these ``bad'' $\beta_k$ values further away from zero. Therefore, we propose our second strategy:
\begin{equation}\label{damp1}
\hat{\beta}_k =\begin{cases}
\beta_k &\text{if $\beta_k \geq \eta$},\\
1-\beta_k &\text{if $\beta_k< \eta$}.\\
\end{cases}
\end{equation}
We note here that there is some trade-off between stability and speed of convergence. This does not mean that larger $\beta_k$ work better, since larger $\beta_k$ may not speed up the convergence if it is not appropriate. Therefore, damping is good, but over-damping may cause inconsistencies and stagnation. In our numerical experiment, we take $\eta=0.3$ in \eqref{damp1} as an example. The results are in \Cref{fig:f4}. Compared with the results in \Cref{fig:f1} and \Cref{fig:new_f6}, it becomes better. We see that there are almost no inconsistencies and there is faster convergence. We also plot the new $\beta_k$ in \Cref{fig:f4}. 

\begin{figure}[htbp]
  \centering
  \includegraphics[height=3.4in, width=4.5in]{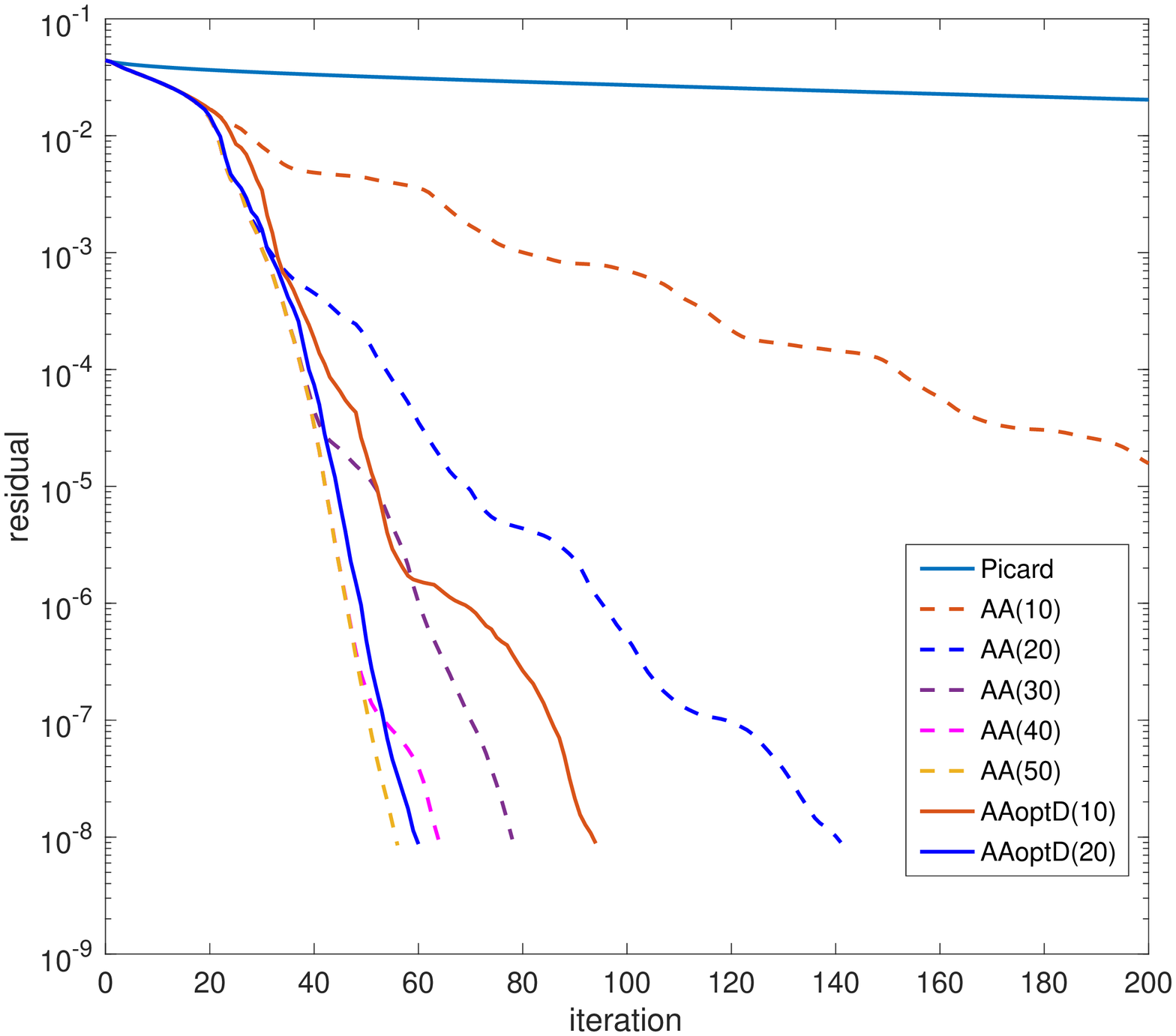}
  \caption{Using larger size windows and bounding the damping factor away from zero.}
  \label{fig:f6}
\end{figure}

\begin{figure}[htbp]
  \centering
  \includegraphics[height=3.4in, width=4.5in]{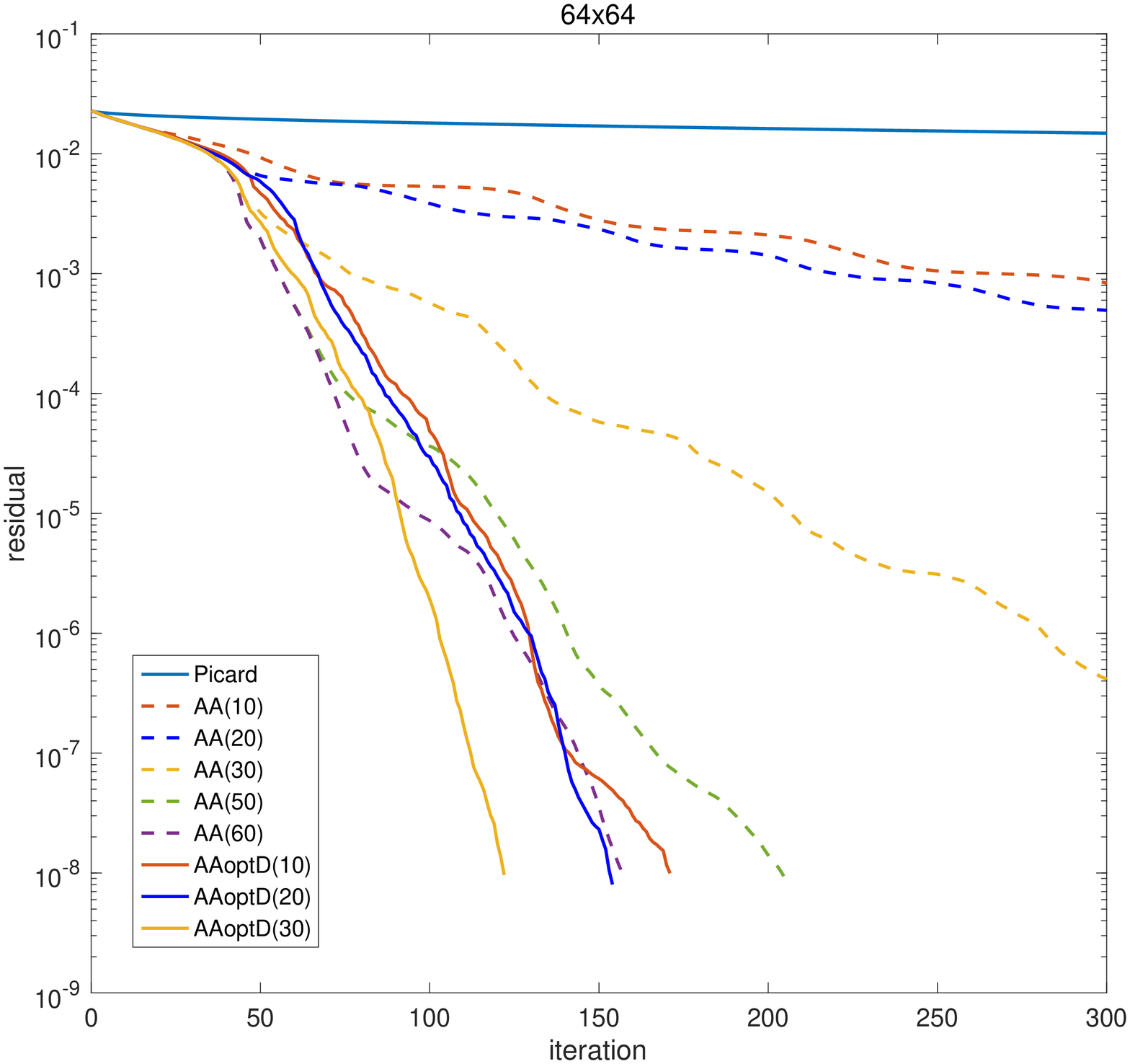}
  \caption{Scaling: solve the Bratu problem on a $64\times 64$ gird.}
  \label{fig:f7_new}
\end{figure}

\begin{figure}[htbp]
  \centering
  \includegraphics[height=3.4in, width=4.5in]{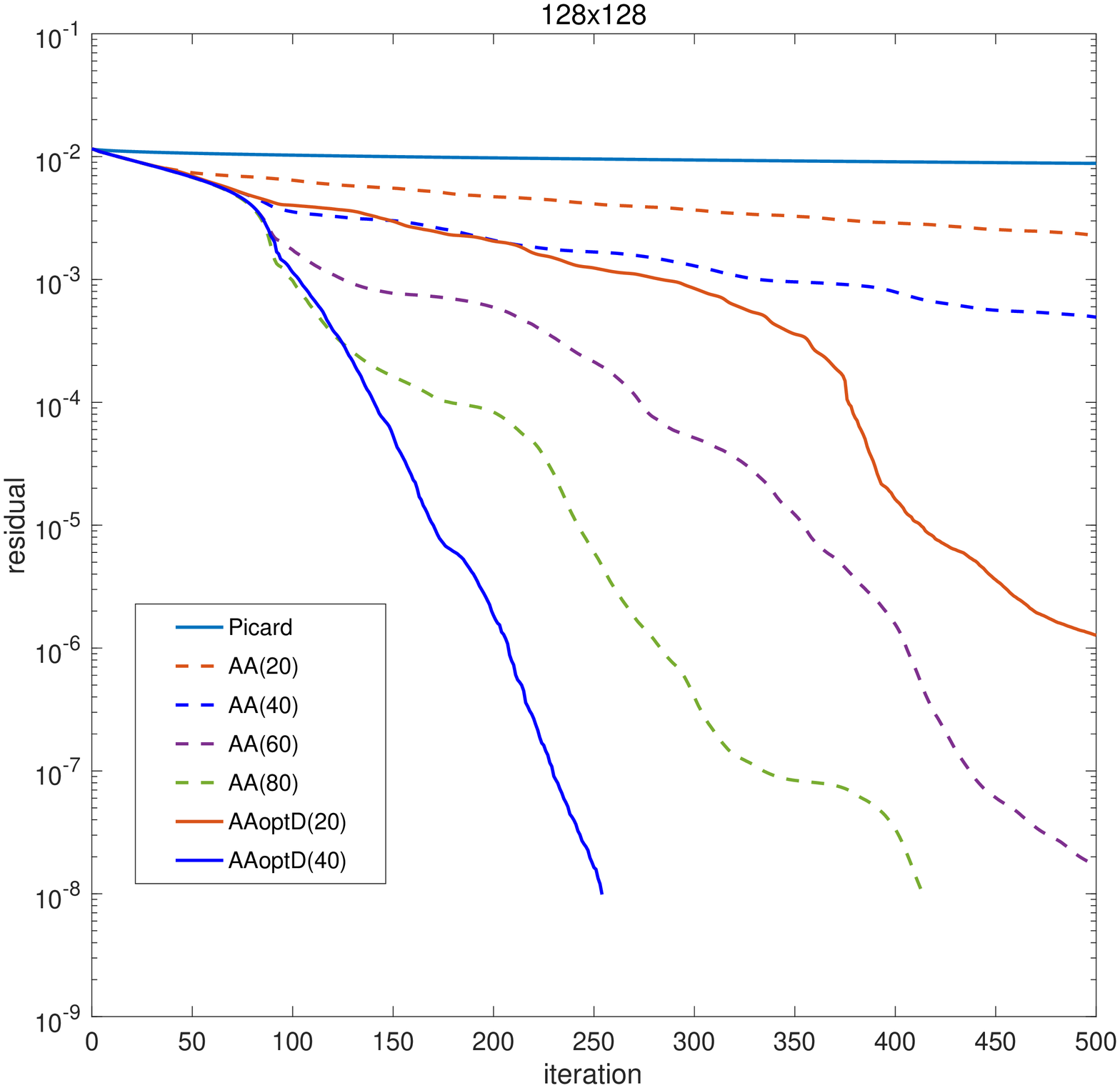}
  \caption{Scaling: solve the Bratu problem on a $128\times 128$ grid.}
  \label{fig:f8_new}
\end{figure}

To compare with the results provided in \cite{WaNi2011}, we go further to increase the windows until $m=50$. Again, without bounding away from zero, there are some stagnations and inconsistencies. To avoid strong over-damping, we apply \eqref{damp1} again with $\eta=0.3$ and obtain our new results in \Cref{fig:f6}. We easily see that $AAoptD(20)$ works as well as $AA(50)$. Moreover, to test its scaling properties, we also solve the Bratu problem on larger grids. In \Cref{fig:f7_new}, for a grid size $64\times 64$, we see that $AAoptD(10)$ is already comparable with $AA(60)$ and $AAoptD(30)$ performs better than $AA(60)$. Similarly, for a grid size $128\times 128$, \Cref{fig:f8_new} shows that $AAoptD(40)$ performs much better than $AA(80)$.

\begin{figure}[htbp]
  \centering
  \includegraphics[height=1.8in, width=4.5in]{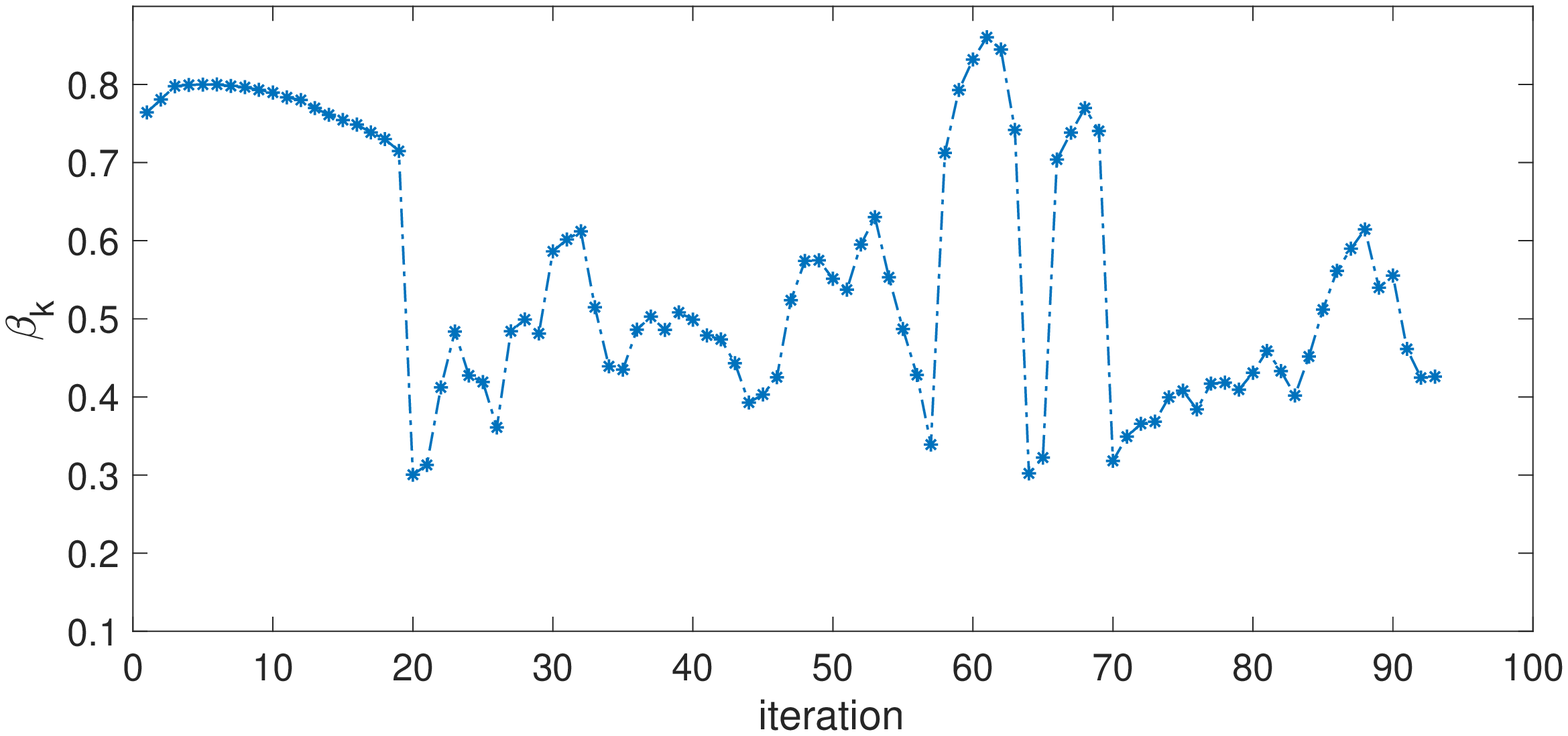}
  \caption{Modified optimal damping factors: $\hat{\beta}_k=1-\beta_k$ when $\beta_k<0.3$.}
  \label{fig:f3}
\end{figure}

\begin{example}\textbf{The\ nonlinear\ convection-diffusion\ problem.} Use AA and AAoptD to solve the following 2D nonlinear convection-diffusion equation in a square region:
$$(-u_{xx}-u_{yy})+(u_x+u_y)+ku^2=f(x,y),\ \  (x,y)\in D=[0,1]\times[0,1] $$
with the source term 
$$f(x,y)=2\pi^2\sin(\pi x)\sin(\pi y)$$
and zero boundary conditions: $u(x,y)=0$ on $\partial D$.
\end{example}

\begin{figure}[htbp]
  \centering
  \includegraphics[height=3.4in, width=4.5in]{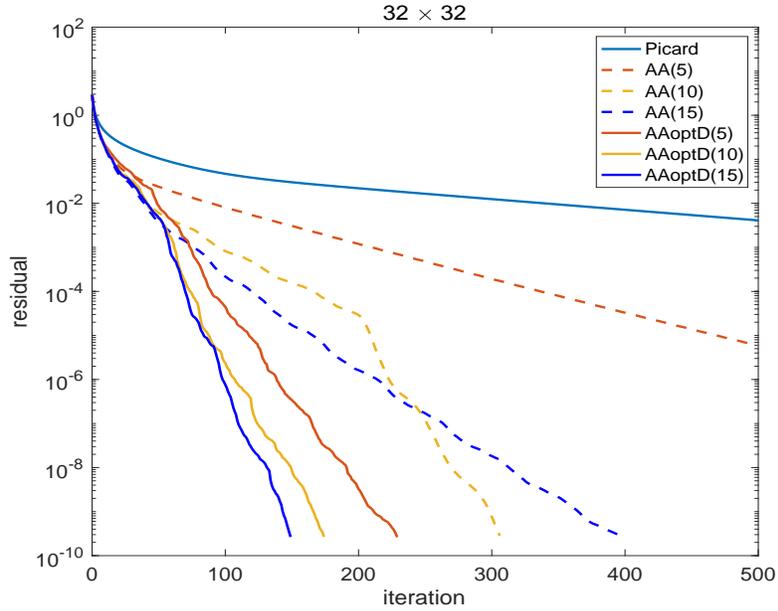}
  \caption{Solving the nonlinear convection-diffusion problem on a $32\times 32$ gird.}
  \label{fig:f11_new}
\end{figure}

\begin{figure}[htbp]
  \centering
  \includegraphics[height=3.4in, width=4.5in]{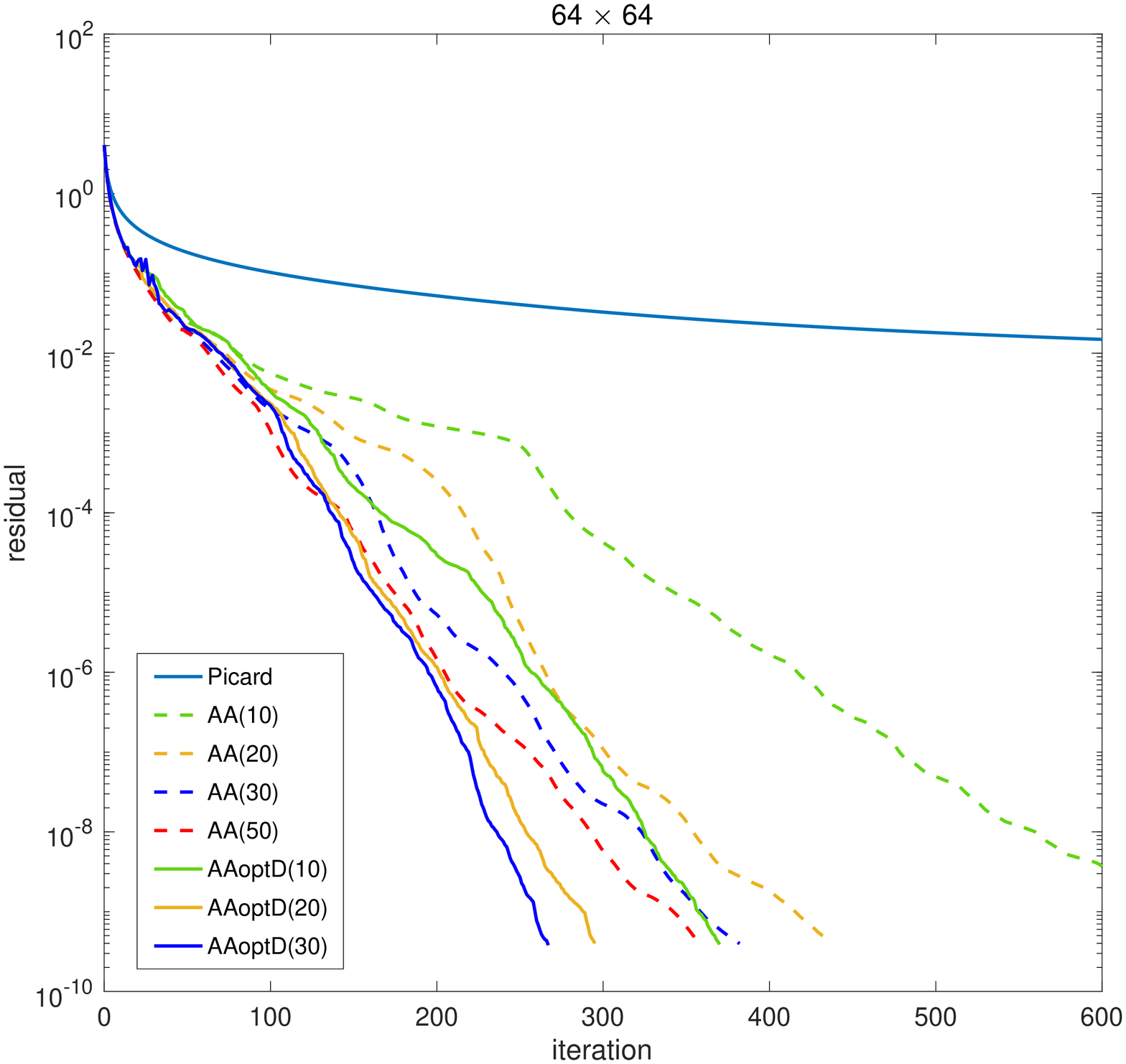}
  \caption{Solving the nonlinear convection-diffusion problem on a $64\times 64$ gird.}
  \label{fig:f12_new}
\end{figure}

In this numerical experiment, we use a centered-difference discretization on $32\times 32$ and $64\times 64$ grids, respectively. We take $k=3$ in the above problem and use $u_0=(1,1,\cdots,1)^{T}$ as an initial approximate solution in all cases. As in solving the Bratu problem, the same preconditioning strategy is used here so that the basic Picard iteration still works. To bound $\beta_k$ away from zero, we use \eqref{damp1} with $\eta=0.25$. The results are shown in \Cref{fig:f11_new} and \Cref{fig:f12_new} for $n=32\times 32$ and $n=64\times 64$, respectively. From \Cref{fig:f11_new}, we see that $AAoptD(5)$ is already better than $AA(15)$; From \Cref{fig:f12_new}, we also observe that $AAoptD(20)$ is better than $AA(50)$. In both cases, $AAoptD(m)$ does a much better job than $AA(m)$, which is consistent with our previous example.

Our next example is about solving a linear system $Ax=b$. As proved by Walker and Ni in \cite{WaNi2011}, AA without truncation is ``essentially equivalent'' in a certain sense to the GMRES method for linear problems.
\begin{example} \textbf{The linear equations.} Apply AA and AAoptD to solve the following linear system $Ax=b$, where $A$ is
\begin{equation*}
A = 
\begin{pmatrix}
2 & -1 & \cdots & 0 & 0 \\
-1 & 2 & \cdots & 0 & 0 \\
\vdots  & \vdots  & \ddots & \vdots & \vdots \\
0 & 0 & \cdots & 2 &-1\\
0 & 0 & \cdots & -1 &2
\end{pmatrix}, \ \ A\in R_{n\times n}
\end{equation*}
and
$$b=(1,\cdots,1)^{T}.$$

Choose $n=10$ and $n=100$, respectively. Here, we choose a large $n$ so that a large window size $m$ is needed in Anderson Acceleration. We also note that the Picard iteration does not work for this problem.
\end{example}

\begin{figure}[htbp]
  \centering
  \includegraphics[height=3.4in, width=4.5in]{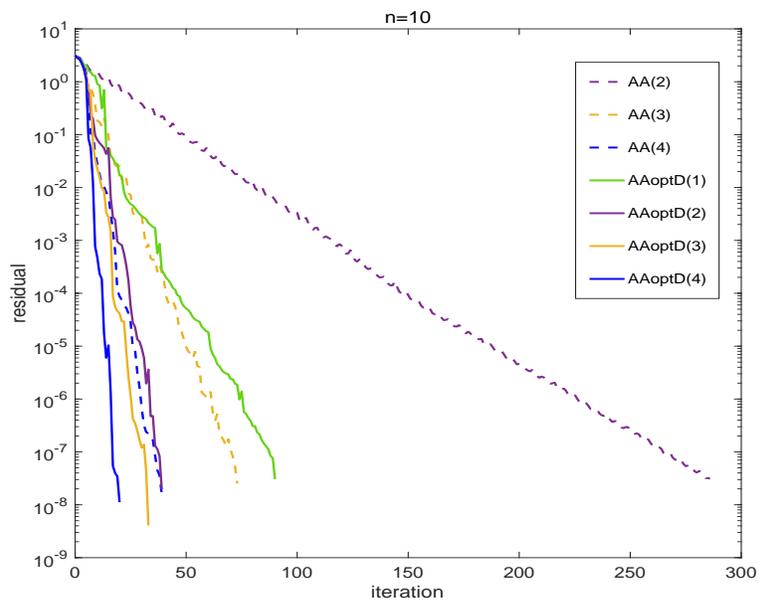}
  \caption{Small $m$: solving a linear problem $Ax=b$ with $n=10$.}
  \label{fig:f13_new}
\end{figure}

\begin{figure}[htbp]
  \centering
  \includegraphics[height=3.4in, width=4.5in]{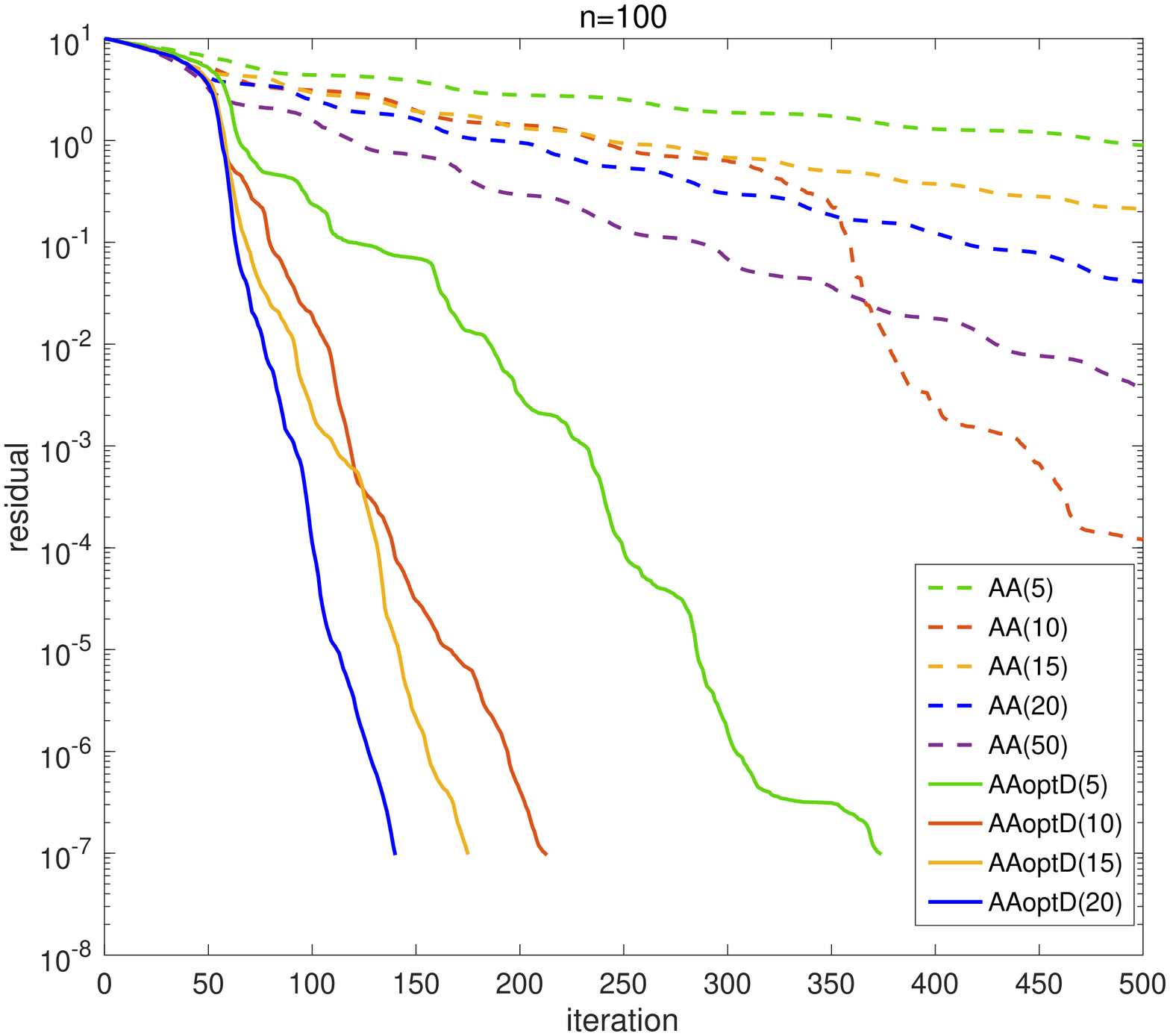}
  \caption{Large $m$: solving a linear problem $Ax=b$ with $n=100$.}
  \label{fig:f14_new}
\end{figure}

The initial guess is $x_0=(0,\cdots,0)^{T}$. Without bounding $\beta_k$ away from zero, the results are shown in \Cref{fig:f13_new} and \Cref{fig:f14_new}. For small $m$, $AA(1)$ does not work, but $AAoptD(1)$ works. Moreover, we obtain from \Cref{fig:f13_new} that $AAoptD(m)$ still does better than $AA(m)$. When $n=100$, we need larger $m$ values. In this case, as shown in \Cref{fig:f14_new}, $AAoptD(5)$ already performs much better than $AA(25)$. This example shows that $AAoptD$ can also be used to solve linear problems.

Finally, we consider cases where very small $m$ works. Our example is from Toth and Kelley's paper \cite{ToKe2015}, where AA is applied to solve the Chandrasekhar H-equation. 
\begin{example} the Chandrasekhar H-equation, arising in Radiative Heat Transfer theory, is a nonlinear integral equation:
$$H(\mu)=\mathbf{G}(H)=\left(1-\frac{c}{2}\int_{0}^{1}\frac{\mu}{\mu+v}H(v)dv\right)^{-1},$$ 
where $c\in [0,1)$ is a physical parameter. 

We will discretize the equation with the composite midpoint rule. Here we approximate integrals on $[0,1]$ by
$$\int_{0}^{1}f(\mu)d\mu\approx\frac{1}{N}\sum_{j=1}^{N}f(\mu_j)$$
where $\mu_j=(i-1/2)/N$ for $1\leq i\leq N$. The resulting discrete problem is 
$$F(x)_i=x_i-\left(1-\frac{c}{2N}\sum_{j=1}^{N}\frac{\mu_ix_j}{\mu_i+\mu_j}\right)^{-1},$$
which is a fully nonlinear system.

It is known \cite{Ke1978} both for the continuous problem and its following midpoint rule discretization, that if $c<1$
$$\rho(\mathbf{G}'(H^*))\leq1-\sqrt{1-c}<1,$$
where $\rho$ denotes spectral radius. Hence the local convergence theory and Picard iteration works.
\end{example}

\begin{figure}[htp]
\centering 
\includegraphics[height=3.4in, width=4.5in]{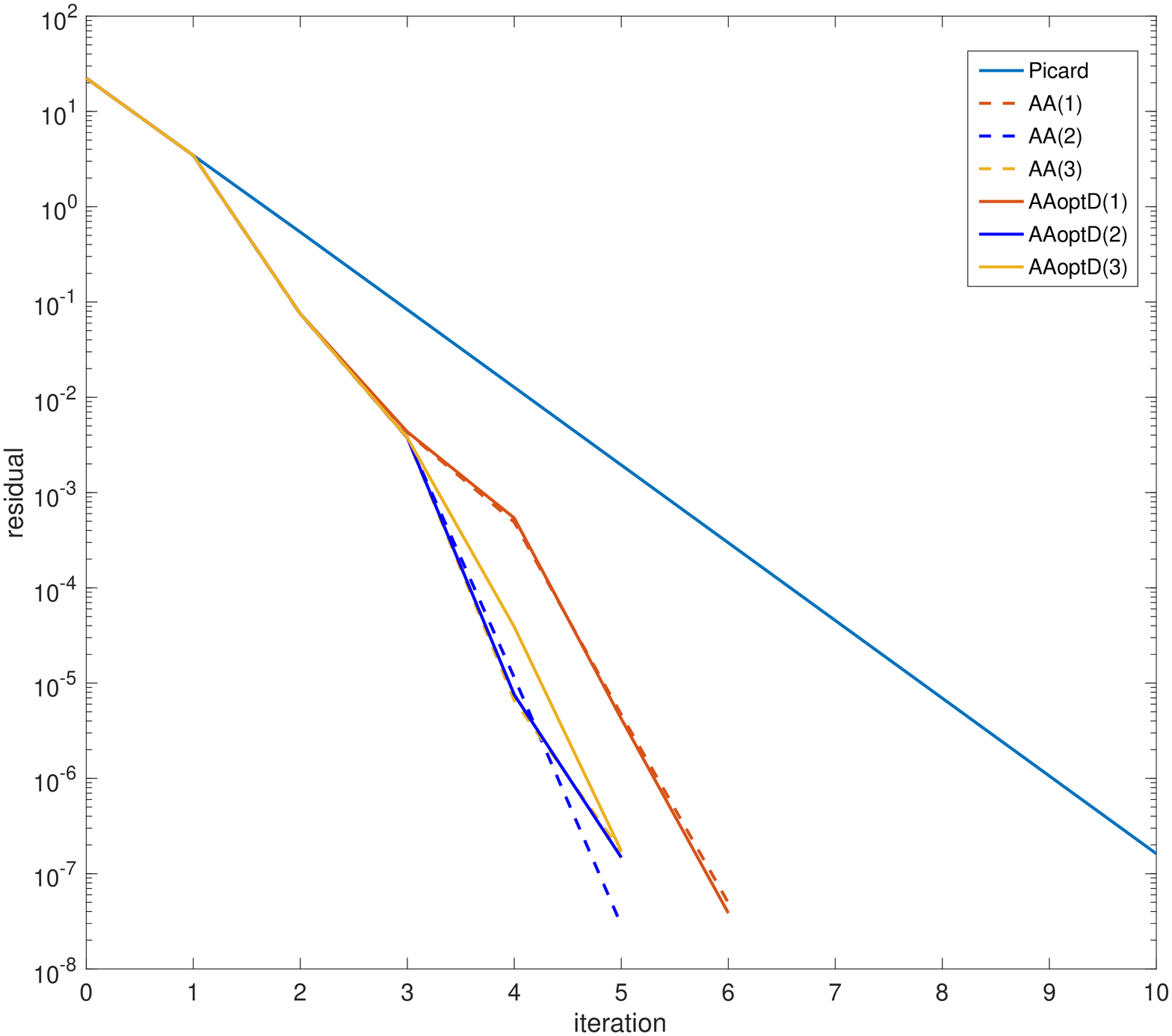}
\caption{Solving Chandrasekhar H-equation with AA and AAoptD: $c=0.5$}
\label{fig:h1}
\end{figure}

\begin{figure}[htp]
\centering 
\includegraphics[height=3.4in, width=4.5in]{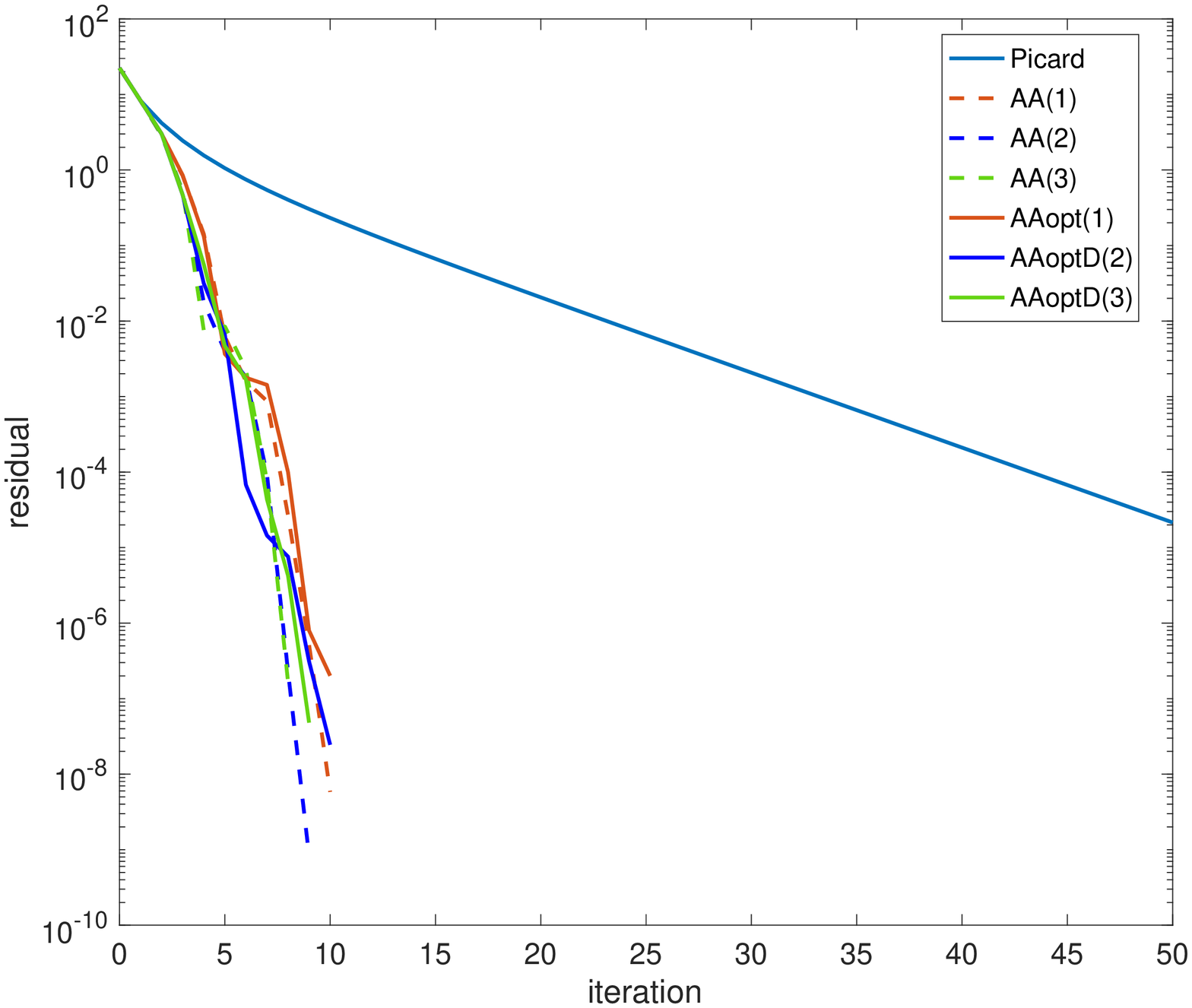}
\caption{Solving Chandrasekhar H-equation with AA and AAoptD: $c=0.99$}
\label{fig:h2}
\end{figure}

\begin{figure}[htp]
\centering 
\includegraphics[height=3.4in, width=4.5in]{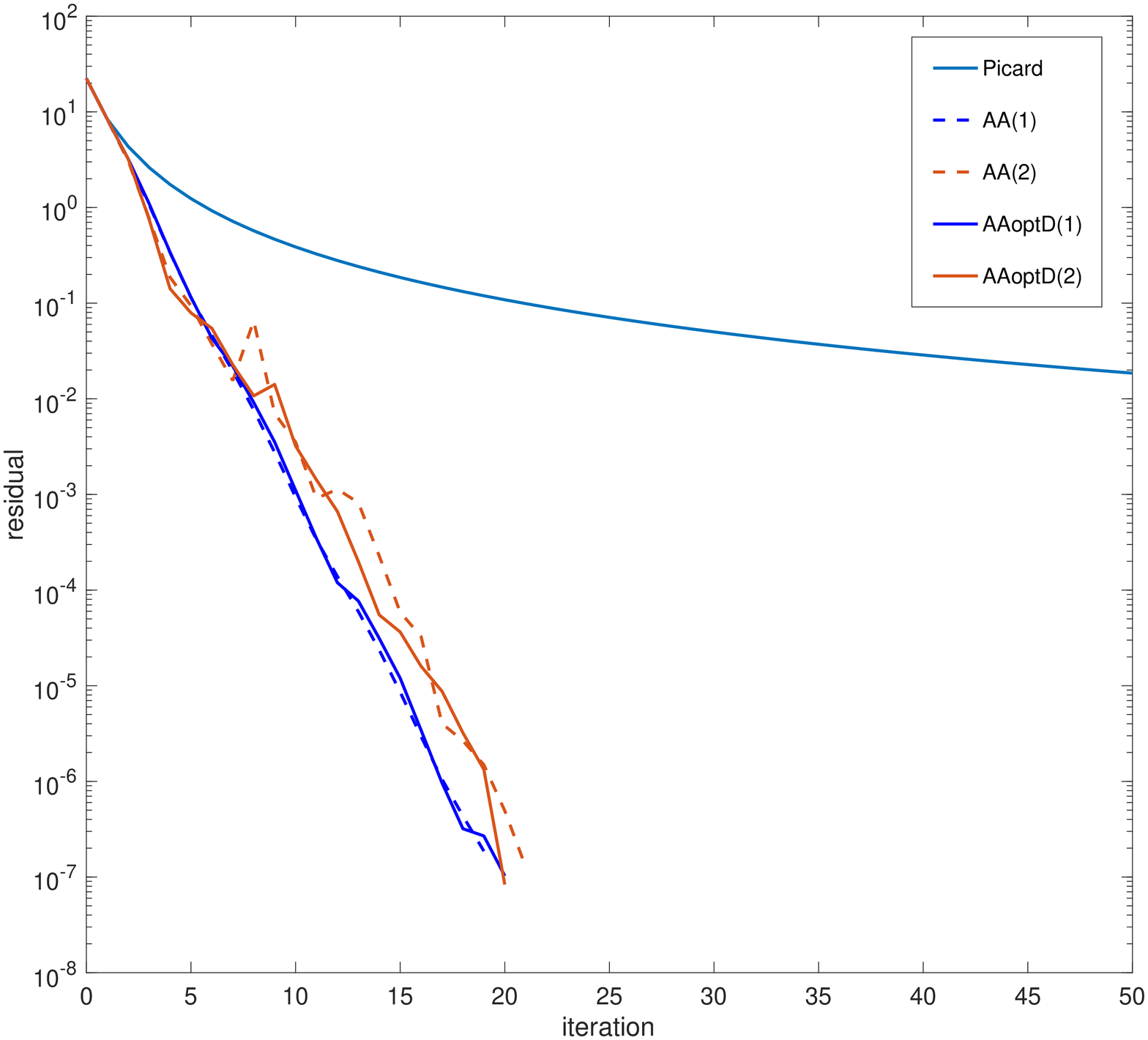}
\caption{Solving Chandrasekhar H-equation with AA and AAoptD: $c=1$}
\label{fig:h3}
\end{figure}

In our numerical experiment, we choose $N=500$, $c=0.5$, $c=0.99$ and $c=1$. The case $c=1$ is a critical value (Picard does not work in this case, but AA does). The numerical results are in \Cref{fig:h1} to \Cref{fig:h3}. Firstly, $AA(m)$ and $AAoptD(m)$, with very small m($\leq 3$) values, work for all cases including the critical case $c=1$ and their performances are comparable. Secondly, increasing $m$ does not always increase the performance. Thirdly, AAoptD may not always have advantages over AA for small window size $m$. This result is reasonable since AAoptD(m) is kind of like packaging $AA(m)$ and $AA(1)$. If $m$ is small, there is almost no difference between $AA(m)$ and $AA(1)$, thus packaging them (varying window sizes) may not give better results.


\section{Conclusions}
\label{sec:conclusions}
We proposed a non-stationary Anderson acceleration algorithm with an optimized damping factor in each iteration to further speed up linear and nonlinear iterations by applying one extra optimization.  This procedure has a strong connection to another perspective of generating non-stationary AA (i.e. varying the window size $m$ at different iterations). It turns out that choosing optimal $\beta_k$ is somewhat similar to packaging sAA(m) and sAA(1) within a single iteration in a cheap way. Moreover, by taking benefit of the QR decomposition in the first optimization problem, the calculation of optimized $\beta_k$ at each iteration is cheap if two extra function evaluations are relatively inexpensive. Our numerical results show that the gain of doing this extra optimized step on $\beta_k$ could be large. Moreover, damping is good but over damping is not good because it may slow down the convergence rate. Therefore, when the stationary AA is not working well or a larger size of the window is needed in AA, we recommend to use AAoptD proposed in the present work.

\section*{Acknowledgments}
This work was partially supported by the National Natural Science Foundation of China [grant number 12001287]; the Startup Foundation for Introducing Talent of Nanjing University of Information Science and Technology [grant number 2019r106]; The first author Kewang Chen also gratefully acknowledge the financial support for his doctoral study provided by the China Scholarship Council (No. 202008320191).

\bibliography{mybibfile}

\end{document}